# Алгоритмические модели поведения человека и стохастическая оптимизация[1]

Бектемесов М.А., д.ф.-м.н., профессор, декан мехмата, Казахский национальный университет имени аль-Фараби,
г. Алматы, Республика Казахстан, E-mail: maktagali@gmail.com
Гасников А.В., д.ф.-м.н., доцент, в.н.с. ИППИ РАН сектор 7, Московский физико-технический институт,
г. Москва, Российская Федерация, E-mail: gasnikov@yandex.ru
Лагуновская А.А., инженер, Московский физико-технический институт,
г. Москва, Российская Федерация,
E-mail: a.lagunovskaya@phystech.edu
Ордабаева Ж.М., студ. Мехмата, КазНУ Казахский национальный университет имени аль-Фараби,
г. Алматы, Республика Казахстан, E-mail: ordabayevazhanna@gmail.com

В статье исследуется параллелизация вычислений при решении задач стохастической оптимизации; рассматривается приложение полученных здесь результатов к поиску равновесного распределения потоков по путям; исследуется зависимость скорости сходимости оптимальных алгоритмов в задачах стохастической безградиентной оптимизации, в зависимости от числа обращений к оракулу за реализацией функции на каждой итерации. Отличительная особенность данной статьи – демонстрация полученных результатов наглядными примерами.

## Адам мінез-құлық және стохастикалық оңтайландырудың алгоритмдік моделдері

Бектемесов М.А., Әл-Фараби атындағы ҚазҰУ механика-математика факультетінің деканы, ф.-м.ғ.д., профессор,
Алматы қ., Қазақстан Республикасы, E-mail: maktagali@gmail.com
Гасников А.В., д.ф.-м.н., доцент, в.н.с. ИППИ РАН 7 сектор, Мәскеу физико-техникалық институты, Мәскеу қ., Ресей Федерациясы,
E-mail: gasnikov@yandex.ru
Лагуновская А.А., инженер, Мәскеу физико-техникалық институты,

[1] Данная статья представляет собой расширенную версию записи доклада, сделанного А.В. Гасниковым на научном семинаре мехмата КазНУ 24 февраля 2017 года.

Мәскеу қ., Ресей Федерациясы, E-mail: a.lagunovskaya@phystech.edu
Ордабаева Ж.М., Мехмат студенті, Әл-Фараби атындағы ҚазҰУ, Алматы қ., Қазақстан Республикасы, E-mail: ordabayevazhanna@gmail.com

Мақалада параллелизация есептеу кезіндегі стохастикалық оңтайландыру зерттеледі; қосымша алынған нәтижелерін іздеуге тепе-теңдік бөлу ағындар жолдары бойынша қаралады; тәуелділік жылдамдығы жинақтылық оңтайлы алгоритмдер міндеттерінде стохастикалық градиентсіз оңтайландыру байланысты өтініштердің санының оракул бақылау функцияларын әрбір итерациясы зерттеледі. Бұл мақаланың ерекшелігі – алынған нәтижелердің көрнекі мысалдар демонстрациясы.

## Algorithmic models of human behavior and stochastic optimization

M.A. Bektemessov., Dr of phys. and math. sc., Dean of the Faculty of Mathematics and Mechanics, Al-Farabi Kazakh National University, Almaty, The Republic of Kazakhstan, E-mail: maktagali@gmail.com
A.V. Gasnikov, Dr of phys. and math. sc., Associated professor of Mathematical Foundation of Control Moscow Institute of Physics and Technology, Moscow, Russia. E-mail: gasnikov@yandex.ru
A.A.Lagunskaya, engineer, Moscow Institute of Physics and Technology, Moscow, Russia, E-mail: a.lagunovskaya@phystech.edu
Zh.M. Ordabayeva, student, Al-Farabi Kazakh National University, Almaty, The Republic of Kazakhstan, E-mail: ordabayevazhanna@gmail.com

The article explores the parallelization of computations in solving stochastic optimization problems; The application of the results obtained here to the search for an equilibrium distribution of flows along paths is considered; The dependence of the rate of convergence of optimal algorithms in the problems of stochastic, gradient optimization is investigated, depending on the number of calls to the oracle behind the implementation of the function at each iteration. A distinctive feature of this article is the demonstration of the results obtained with illustrative examples.

## 1. Введение

Настоящая статья носит обзорный характер. Ее цель продемонстрировать несколькими яркими (на наш взгляд) примерами то, как можно интерпретировать ряд современных результатов о сходимости (стохастического) метода зеркального спуска и его безградиентных вариантов.

В разделе 2 описывается стохастический зеркальный спуск – численный метод решения задач выпуклой стохастической оптимизации. Исследуется вопрос о возможности его параллелизации.

В разделе 3 результаты раздела 2 переносятся на безградиентные методы – нет возможности считать (стохастический) субградиент.

В разделе 4 результаты раздела 2 переносятся на конкретную задачу – поиск равновесного распределения транспортных потоков.

## 2. Стохастический зеркальный спуск

В данном разделе описывается метод зеркального спуска – МЗС (А.С. Немировский, 1977 [1, 2]). Обсуждается стохастический и онлайн варианты метода, а также возможность его параллелизации.

Рассматривается класс задач выпуклой стохастической оптимизации

$$f(x) = E_\xi\left[f(x,\xi)\right] \to \min_{x\in Q}, \tag{1}$$

где $f(x)$ – выпуклая функция, а $Q$ – выпуклое множество. Часто в приложениях, рассматриваемых в данной статье

$$Q = S_n(1) = \left\{x\in\mathbb{R}_+^n:\ \sum_{i=1}^n x_i = 1\right\}.$$

На каждой итерации (шаге) можно посчитать в выбранной точке $x^k$ стохастический субградиент $\partial_x f\left(x^k,\xi^k\right)$, где $\left\{\xi^k\right\}_k$ – независимые одинаково распределенные случайные величины (распределенные также как $\xi$). Для простоты рассуждений будем считать, что работаем с каким-то конкретным измеримым селектором $\nabla_x f\left(x^k,\xi^k\right)$. При этом

$$E_\xi\left[\nabla_x f(x,\xi)\right] \equiv \nabla f(x). \tag{2}$$

Вводится $p$-норма ($p\in[1,2]$) и параметр $q$: $1/p + 1/q = 1$. Предполагается, что существует такое $M>0$, что

$$E_\xi\left[\left\|\nabla_x f(x,\xi)\right\|_q^2\right] \le M^2,\ q\in[2,\infty]. \tag{3}$$

Вводится прокс-функция $d(x) \geq 0$ ($d(x^0) = 0$, в ряде случаев точку старта будем обозначать не $x^0$, а $x^1$, тогда $d(x^1) = 0$), являющаяся 1-строго выпуклой функцией в $p$-норме. Вводится расхождение Брэгмана

$$V(x,y) = d(x) - d(y) - \langle \nabla d(y), x - y \rangle.$$

Популярными примерами таких функций являются [2]:

$$p = 1,\ d(x) = \ln n + \sum_{k=1}^{n} x_k \ln x_k,\ V(x,y) = \sum_{k=1}^{n} x_k \ln(x_k / y_k),\ Q = S_n(1),$$

$$p = 2,\ d(x) = \frac{1}{2}\|x\|_2^2,\ V(x,y) = \frac{1}{2}\|x - y\|_2^2,\ Q = \mathbb{R}^n.$$

Далее вводится оператор проектирования согласно расхождению Брэгмана

$$\mathrm{Mirr}_{x^k}(\mathrm{v}) = \arg\min_{x \in Q} \left\{ \langle \mathrm{v}, x - x^k \rangle + V(x, x^k) \right\}.$$

С помощью этого оператора вводится МЗС

$$\boxed{x^{k+1} = \mathrm{Mirr}_{x^k}\left(h \nabla_x f(x^k, \xi^k)\right),\ k = 0, ..., N-1},$$

где

$$h = \frac{R}{M}\sqrt{\frac{2}{N}} = \frac{\varepsilon}{M^2},\ R^2 = V(x_*, x^0),$$

здесь $x_*$ – решение задачи (1). Если $x_*$ не единственно, то можно считать, что среди всех решений (1) выбирается такое $x_*$, которое доставляет минимум $V(x_*, x^0)$.

**Теорема 1 (А.С. Немировский [2]).** *Пусть выполняются условия (2), (3). Тогда после*

$$N = \frac{2M^2 R^2}{\varepsilon^2}$$

*итераций (вычислений стохастических градиентов)*

$$E\left[ f(\bar{x}^N) \right] - f_* \leq \varepsilon,$$

*где*

$$f_* = \min_{x \in Q} f(x),\ \bar{x}^N = \frac{1}{N}\sum_{k=0}^{N-1} x^k.$$

Следующая теорема немного уточняет в случае тяжелых хвостов стохастического субградиента (см. условие iii)) результат Немировского–Юдицкого–Лана–Шапиро, 2009 [3].

**Теорема 2.** *Пусть выполняется условие (2) и одно из следующих условий*

*i)* $\left\|\nabla_x f(x,\xi)\right\|_q \le M$ *для всех* $x \in Q$ *и п.н. по* $\xi$;

*ii)* $E_\xi\left[\exp\left(\left\|\nabla_x f(x,\xi)\right\|_q^2 \big/ M^2\right)\right] \le \exp(1)$, $\ln\sigma^{-1} \ll N$;

*iii)* *Существует такое* $\alpha > 2$, *что*

$$P\left(\left\|\nabla_x f(x,\xi)\right\|_q^2 \big/ M^2 \ge t\right) \le (1+t)^{-\alpha},\ \sigma^{-1/(\alpha-1)} \ll N.$$

*Тогда при соответствующих* $\sigma$ *(см. ii) и iii)) имеет место неравенство*

$$P\left( f\left(\bar{x}^N\right) - f_* \le \frac{C_1 M}{\sqrt{N}}\left(R + C_2 \bar{R}\sqrt{\ln\sigma^{-1}}\right)\right) \ge 1-\sigma,$$

*где* $\bar{R} = \max_{x,y\in Q}\|x-y\|_p$, $C_1, C_2 \le 3$ *за исключением случая iii), в котором* $C_1(\alpha)$.

С помощью теоремы 2 был получен следующий результат [4], который может быть использован для распараллеливания работы МЗС.

**Теорема 3.** *Пусть выполняется условие (2) и одно из условий i) – iii). Тогда существует такое* $C \le 50C_1$, *что по* $K = \lceil 2\log_2\sigma^{-1}\rceil$ *независимым друг от друга в смысле выбора реализаций* $\left\{\xi^k\right\}_k$ *траекториям МЗС (на каждой траектории рассчитываем* $\bar{x}^{N,i}$*) с одинаковым числом шагов*

$$N = \frac{CM^2\bar{R}^2}{\varepsilon^2}$$

*можно "собрать" такой*

$$\bar{x}^N = \frac{1}{K}\sum_{i=1}^{K}\bar{x}^{N,i},$$

*что*

$$P\left( f\left(\bar{x}^N\right) - f_* \le \varepsilon\right) \ge 1-\sigma.$$

**Пример 1.** Результат теоремы 3 допускает следующую интерпретацию. Пусть есть один мудрец, который прожил 1000 лет и есть 10 экспертов, проживших по 100 лет. Считается, что у всех у них одинаковая функция ежедневных затрат $f(x)$ и каждый стремится оказаться как можно быстрее в минимуме этой функции. Каждый день эксперт или мудрец (оптимально) выбирают способ “как жить”, исходя из своей истории. Оказывается, что при весьма общих условиях на “обратную связь”, которую можно получать ежедневно ($\nabla_x f(x^k,\xi^k)$), мудрец “эквивалентен” соответствующему числу независимых экспертов в смысле итогового понимания, как именно надо было жить. Под эквивалентностью понимается “с точностью до универсального числового множителя”, и результат, с небольшими оговорками, не зависит от того, сколько именно прожил мудрец и эксперты. Важно только, чтобы эксперты в сумме прожили приблизительно столько же сколько мудрец. При переходе на такую интерпретацию в теореме 3, по сути, получено решение задачи, поставленной Ю.Е. Нестеровым в 2003 году [5]. □

Далее приведенные выше результаты распространяются, следуя [6], на задачи, в которых дополнительно известно, что

$$f(x) - \mu_2\text{-сильной выпуклая функция в 2-норме}, \tag{4}$$

при этом выбирается $p = 2$. В варианте аналогичном теореме 1 это также было ранее сделано А.С. Немировским [1]. Более того, в 2011 году А.С. Немировским и А.Б. Юдицким теорема 1 была обобщена на случай, когда вместо условия (2) выполняется более слабое условие [7]

$$\left\| E_\xi\left[\nabla_x f(x,\xi^k)\right] - \nabla f(x) \right\|_q \le \delta/\bar{R}. \tag{5}$$

**Теорема 4 (Немировский–Юдицкий).** *Пусть выполняются условия (3), (5). Тогда*

$$E\left[f(\bar{x}^N)\right] - f_* \le MR\sqrt{\frac{2}{N}} + \delta,$$

*где*

$$\bar{x}^N = \frac{1}{N}\sum_{k=0}^{N-1} x^k.$$

*Пусть выполняются условия (3), (4), (5) ($p = 2$). Тогда*

$$E\left[f\left(\bar{x}^N\right)\right]-f_*\le\frac{M_2^2}{2\mu_2 N}\left(1+\ln N\right)+\delta$$

*где*

$$\bar{x}^N=\frac{1}{N}\sum_{k=1}^{N}x^k,\ x^{k+1}=\operatorname{Mirr}_{x^k}\left(h_k\nabla_x f\left(x^k,\xi^k\right)\right),\ h_k=\left(\mu_2 k\right)^{-1},\ k=1,\dots,N.$$

В работе [8] удалось еще больше ослабить условие (5)

$$\sup_{\left\{x^k=x^k\left(\xi^1,\dots,\xi^{k-1}\right)\right\}_{k=1}^{N}\in Q} E\left[\frac{1}{N}\sum_{k=1}^{N}\left\langle E_{\xi^k}\left[\nabla_x f\left(x^k,\xi^k\right)\middle|\xi^1,\dots,\xi^{k-1}\right]-\nabla f\left(x^k\right),x^k-x_*\right\rangle\right]\le\delta. \quad (6)$$

При этом формулировка теоремы 4 полностью сохранилась. Именно в таком варианте с наиболее общим условием (6) на уровень допустимого шума $\delta$ в следующем разделе удалось воспользоваться теоремой 4 для получения оценок скорости сходимости безградиентных методов.

Приведенные выше результаты можно обобщить на задачи онлайн оптимизации [9]. Оказывается, что теорема 4 имеет место с точно такими же оценками скорости сходимости и для задач онлайн оптимизации. В условиях шума (6) ($\delta>0$) такие результаты были получены впервые.

Оценки теоремы 4 являются неулучшаемыми в выпуклом случае с точностью до мультипликативной константы и с точностью до множителя вида $\sim\ln N$ в сильно выпуклом случае. В случае задач онлайн оптимизации ($\delta=0$) оценки теоремы 4 являются в обоих случаях неулучшаемыми с точностью до мультипликативной константы [10].

Приведем вкратце постановку задачи онлайн оптимизации и полученный результат. Положим (все функции $f_k\left(x\right)$ из того же класса, что была функция $f\left(x\right)$, причем $f_k\left(x\right)$ может враждебно подбираться, исходя из знания истории $\left\{x^1,\dots,x^{k-1}\right\}$)

$$\operatorname{Regret}_N\left(\left\{f_k\left(\cdot\right)\right\},\left\{x^k\right\}\right)=\frac{1}{N}\sum_{k=1}^{N}f_k\left(x^k\right)-\min_{x\in Q}\frac{1}{N}\sum_{k=1}^{N}f_k\left(x\right)$$

и рассмотрим следующую задачу

$$\frac{1}{N}\sum_{k=1}^{N}f_k\left(x^k\right)\to\min_{\left\{x^k\left(\nabla_x f_1\left(x^1,\xi^1\right);\dots;\nabla_x f_{k-1}\left(x^{k-1},\xi^{k-1}\right)\right)\in Q\right\}_{k=1}^{N}}.$$

Если в теореме 4 в процедуре МЗС при расчете $x^{k+1}$ вместо $\nabla_x f\left(x^k,\xi^k\right)$ использовать $\nabla_x f_k\left(x^k,\xi^k\right)$ – то, что доступно наблюдению только на этом шаге (можно называть этот вектор $\delta$-смещенным стохастическим субградиентом $f_k\left(x\right)$), то теорема 4 останется верной если заменить

$$E\left[f\left(\bar{x}^N\right)\right]-f_* \text{ на } E\left[\text{Regret}_N\left(\left\{f_k\left(\cdot\right)\right\},\left\{x^k\right\}\right)\right].$$

Приводимый далее (известный ранее) конкретный пример задачи онлайн оптимизации (взвешивание экспертных мнений) будет играть важную роль в разделе 4.

**Пример 2.** Если $f_k\left(x\right)=\left\langle l^k,x\right\rangle$, $Q=S_n\left(1\right)=\left\{x\in\mathbb{R}^n_+:\ \sum_{i=1}^n x_i=1\right\}$, то стоит выбирать $p=1$, $V\left(x,y\right)=\sum_{k=1}^n x_k\ln\left(x_k/y_k\right)$. В этом случае

$$x_i^{k+1}=\frac{\exp\left(-h\sum_{r=1}^k l_i^r\right)}{\sum_{j=1}^n\exp\left(-h\sum_{r=1}^k l_j^r\right)}=\frac{x_i^k\exp\left(-h\cdot l_i^k\right)}{\sum_{j=1}^n x_j^k\exp\left(-h\cdot l_j^k\right)},$$

$$\text{Regret}_N\left(\left\{f_k\left(\cdot\right)\right\},\left\{x^k\right\}\right)\le MR\sqrt{\frac{2}{N}},$$

где

$$h=\frac{R}{M}\sqrt{\frac{2}{N}},\ M=\max_{k=1,\dots,N}\left\|l^k\right\|_\infty,\ R^2=\ln n.$$

Все приведенные выше результаты (в том числе и этот) можно обобщить на случай, когда $h$ выбирается адаптивным, т.е. не зависящим от $N$, образом. Для этого нужно использовать вместо МЗС метод двойственных усреднений (Ю.Е. Нестеров, 2009 [11]). Для большей наглядности мы ограничились здесь изложением варианта с МЗС.

Приведенная оценка регрета также является неулучшаемой с точностью до мультипликативной константы [9].

Полученная выше оценка регрета останется верной (в среднем), если на каждом шаге $k$ в оценке регрета выбирать не $x^k$, а независимо и случайно осуществлять выбор одной из вершин симплекса, согласно распределению вероятностей $\left\{x_i^k\right\}_{i=1}^n$, определенному выше, и использовать то, что получится в качестве $\left\{x^k\right\}_k$, формирующих регрет.

Этот результат можно проинтерпретировать, например, следующим образом. Пусть мы играем с казино в орлянку: каждый день казино сообщает нам один из двух возможных исходов $\{1,-1\}$. Каждый день мы (заранее) должны делать ставку 1 рубль на один из исходов. Если мы угадали, то выигрываем рубль – нам возвращают этот рубль вместе с еще одним, если не угадали, то проигрываем (этот) рубль. Казино не знает нашу текущую ставку, но зато знает все предыдущие наши ставки (и может это использовать враждебным для нас способом). Будем формировать вектор $l^k \in \mathbb{R}^2$ следующим образом: $l^k = (-1,1)^T$ если в $k$-й день казино сообщило 1 и $l^k = (1,-1)^T$ в противном случае. Приведенный выше результат позволяет обыгрывать такое казино, если казино не в состоянии генерировать случайную последовательность $\{1,-1\}$ в следующем смысле: за $N \gg 1$ дней доли (частоты) встречаемости 1 и $-1$ отличаются не более чем на $N^{-1/2}$. □

## 3. Стохастический зеркальный спуск и безградиентная оптимизация

В данном разделе результаты предыдущего раздела переносятся на задачи стохастической онлайн оптимизации, в которых вместо (зашумленного) стохастического градиента доступными наблюдению на каждом шаге (итерации) являются только (зашумленные) реализации значений функции в нескольких точках или зашумленная реализация производной по направлению. Для простоты формулировок далее в этом разделе мы ограничимся только не онлайн случаем и сходимостью в среднем. Основная идея [1, 8] заключается в введении

$$\nabla_x f\left(x^k, \xi^k := \left(\xi^k, e^k\right)\right) := \frac{n}{\tau} f\left(x^k + \tau e^k, \xi^k\right) e^k, \tag{7}$$

(одноточечная обратная связь)

$$\nabla_x f\left(x^k, \xi^k\right) := \frac{n}{\tau}\left( f\left(x^k + \tau e^k, \xi^k\right) - f\left(x^k, \xi^k\right)\right) e^k, \tag{8}$$

(двухточечная обратная связь)

$$\nabla_x f\left(x^k,\xi^k\right) := n\left\langle \nabla_x f\left(x^k,\xi^k\right), e^k \right\rangle e^k . \qquad (9)$$

(обратная связь в виде реализации производной по направлению)

Считается, что $\left\{e^k\right\}_k$ – независимые одинаково распределенные случайные векторы. Будем предполагать, что $f\left(x^k,\xi^k\right)$ в случаях (7), (8) доступны с, вообще говоря, неслучайным шумом уровня $\delta$:

$$\left| E_\xi\left[ f\left(x,\xi\right)\right] - f\left(x\right)\right| \le \delta . \qquad (10)$$

В случае (9) будем считать, что выполняется условие (2). Здесь и далее в этом разделе будем считать, что все условия выполняются не только в множестве $Q$, но и в $\delta$-окрестности (в $2$-норме) этого множества.

В литературе давно открытым стоял вопрос: как следует выбирать $e^k$? Основные конкурирующие способы:

$$e^k \in RS_2^n\left(1\right),$$

т.е. $e^k$ выбирается равновероятно на поверхности евклидовой сферы единичного радиуса в $\mathbb{R}^n$;

$$e^k = (\underbrace{0,\ldots,0,1}_{i},0,\ldots,0) \text{ с вероятностью } 1/n ,$$

покомпонентная рандомизация.

В работе [12] было установлено, что абсолютно во всех случаях, если ориентироваться только на число сделанных итераций, и не смотреть на "стоимость" итерации, первый способ $e^k \in RS_2^n\left(1\right)$ доминирует все известные остальные способы – с точностью до логарифмического множителя, причем в большинстве случаев такая оговорка не требуется.

Заметим, что

$$E_{e^k}\left[ n\left\langle \nabla_x f\left(x^k,\xi^k\right), e^k\right\rangle e^k \right] = \nabla_x f\left(x^k,\xi^k\right), \text{ (следует сравнить с (2))}$$

$$E\left[ \left\| \frac{n}{\tau}\left( f\left(x^k+\tau e^k,\xi^k\right) - f\left(x^k,\xi^k\right)\right) e^k \right\|_q^2 \right] \le \frac{3}{4} n^2 \tau^2 L_2^2 E_{e^k}\left[ \left\| e^k \right\|_q^2 \right] +$$

$$+3n^2 E\left[\left\langle \nabla_x f\left(x^k,\xi^k\right),e^k\right\rangle^2 \left\|e^k\right\|_q^2\right]+12\frac{\delta^2 n^2}{\tau^2}E_{e^k}\left[\left\|e^k\right\|_q^2\right].$$ (следует сравнить с (3))

Если

$$E_{\xi^k}\left[\left|f\left(x,\xi^k\right)\right|^2\right]\le B^2, \tag{11}$$

то

$$E\left[\left\|\frac{n}{\tau}f\left(x^k+\tau e^k,\xi^k\right)e^k\right\|_q^2\right]\le\frac{n^2B^2}{\tau^2}E_{e^k}\left[\left\|e^k\right\|_q^2\right].$$ (следует сравнить с (3))

Также заметим, что мы не можем стремить $\tau\to 0+$ в (7), (8) из-за присутствия в оценках неограниченно возрастающего отношения $\delta/\tau$. Если бы $\delta=0$, то можно было объединить случаи (8) и (9), поскольку при $\tau\to 0+$ (8) переходит в (9).

Используя явление концентрации меры на сфере (П. Леви, А. Пуанкаре, В.Д. Мильман [13, 14]) для $e\in RS_2^n(1)$ можно получить следующие результаты

$$E\left[\left\|e\right\|_q^2\right]\le\min\left\{q-1,16\ln n-8\right\}\cdot n^{\frac{2}{q}-1},\ E\left[\left\langle c,e\right\rangle^2\right]\le\left\|c\right\|_2^2 n^{-1},\ 2\le q\le\infty,$$

$$E\left[\left\langle c,e\right\rangle^2\left\|e\right\|_q^2\right]\le\sqrt{3}\left\|c\right\|_2^2\min\left\{2q-1,32\ln n-8\right\}\cdot n^{\frac{2}{q}-2},\ 2\le q\le\infty.$$

Выписанные формулы позволяют надеяться, что можно погрузить рассматриваемые в этом разделе задачи в класс задач обычной стохастической оптимизации, рассматриваемый в прошлом разделе. И, действительно, уже прямо сейчас это можно сделать для случая (9). Однако для случаев (7), (8) даже при $\delta=0$ возникает проблема с условием (2)

$$E_{e^k}\left[\frac{n}{\tau}f\left(x^k+\tau e^k,\xi^k\right)e^k\right]=E_{e^k}\left[\frac{n}{\tau}\left(f\left(x^k+\tau e^k,\xi^k\right)-f\left(x^k,\xi^k\right)\right)e^k\right]\ne$$

$$\ne E_{e^k}\left[n\left\langle\nabla_x f\left(x^k,\xi^k\right),e^k\right\rangle e^k\right].$$

Выручает введенное ранее условие (6). Благодаря именно такому обобщению теоремы 4 удалось получить следующий результат [15].

**Теорема 5**. *Пусть для задачи выпуклой стохастической оптимизации (1) на каждой итерации доступен наблюдению один из следующих векторов (7) – (9). При этом в случае (7) выполняются условия (3) с* $q=2$*, (10), (11), в случае (8) выполняются условия (3) с* $q=2$*, (10) и условие*

$$\left\|\nabla f(y)-\nabla f(x)\right\|_2 \le L_2\left\|y-x\right\|_2, \quad (12)$$

*в случае (9) выполняются условия (2), (3) с* $q=2$. *Тогда используя МЗС с* $\nabla_x f\left(x^k,\xi^k\right)$, *рассчитываемым по соответствующей формуле (7) – (9), с явно выписываемыми оценками на* $\tau$, $h$ *и допустимым уровнем шума* $\delta \le \delta(\varepsilon,n)$, *после* $N$ *итераций можно получить следующую оценку*

$$E\left[f\left(\bar{x}^N\right)\right]-f_* \le \varepsilon,$$

*где* $N$ *определяется из таблицы 1 в случае (7) и из таблицы 2 в случае (8), (9). Причем,* $\tilde{O}(\ )=O(\ )$ *с точностью до логарифмического по* $n$ *или* $N$ *множителя. Параметр* $p$ *выбирается из отрезка* $[1,2]$, $1/p+1/q=1$.

Таблица 1

| $N$ $\left(R^2=\tilde{O}\left(\left\|x_*-x^0\right\|_p^2\right)\right)$ | $E\left[\left\|\nabla_x f(x,\xi)\right\|_2^2\right]\le M_2^2$ | $\left\|\nabla f(y)-\nabla f(x)\right\|_2 \le L_2\left\|y-x\right\|_2$ |
|---|---|---|
| $f(x)$ выпуклая | $\tilde{O}\left(\frac{B^2M_2^2R^2n^{1+2/q}}{\varepsilon^4}\right)$ | $\tilde{O}\left(\frac{B^2L_2R^2n^{1+2/q}}{\varepsilon^3}\right)$ |
| $f(x)$ – $\mu_2$-сильно выпуклая | $\tilde{O}\left(\frac{B^2M_2^2n^2}{\mu_2\varepsilon^3}\right)$ | $\tilde{O}\left(\frac{B^2L_2n^2}{\mu_2\varepsilon^2}\right)$ |

Таблица 2

| $N$ $\left(R^2=\tilde{O}\left(\left\|x_*-x^0\right\|_p^2\right)\right)$ | $E\left[\left\|\nabla_x f(x,\xi)\right\|_2^2\right]\le M_2^2$ | $\left\|\nabla f(y)-\nabla f(x)\right\|_2 \le L_2\left\|y-x\right\|_2$ |
|---|---|---|
| $f(x)$ выпуклая | $\tilde{O}\left(\frac{M_2^2R^2n^{2/q}}{\varepsilon^2}\right)$ | $\tilde{O}\left(\frac{M_2^2R^2n^{2/q}}{\varepsilon^2}\right)$ |
| $f(x)$ – $\mu_2$-сильно выпуклая | $\tilde{O}\left(\frac{M_2^2n}{\mu_2\varepsilon}\right)$ | $\tilde{O}\left(\frac{M_2^2n}{\mu_2\varepsilon}\right)$ |

**Замечание 1.** В категориях $\tilde{O}(\ )$ даже при $\delta=0$ выписанные в таблицах 1, 2 оценки оптимальны, т.е. совпадают с нижними оценками [16]. В случае таблицы 1 нужны оговорки, что все-таки приведенные оценки могут быть улучшен по тому как входит в оценку $N$ желаемая точность $\varepsilon$ за счет ухудшения того как входит размерность пространства $n$ [17]. В случае таблицы 2 оценки равномерно оптимальны одновременно и по $\varepsilon$ и по $n$. Причем эти оценки остаются оптимальными также в случае, когда рассматривается задача не стохастической, а обычной оптимизации. В точности такие же оценки можно получить и для задач (стохастической) онлайн оптимизации. Иногда в таком контексте, т.е. с обратной связью (7) – (9), эти задачи называют задачами о (контекстуальных) нелинейных многоруких бандитах [18]. Новизна в теореме 5 заключается в рассмотрении случая $\delta > 0$, а также в рассмотрении случаев с $p \neq 2$. Впрочем, недавно в работе [19, 20] были получены результаты, близкие к приведенным в теореме 5 при $p = 1$ и $\delta = 0$.

**Замечание 2.** В случае обратной связи (8) условие на уровень шума $\delta$ имеет вид

$$\delta \leq \min\left\{\frac{\varepsilon\tau}{16R\sqrt{n}}, \frac{M_2\tau}{\sqrt{96n}}\right\},$$

Откуда из условия на $\tau$

$$\tau = \min\left\{\max\left\{\frac{\varepsilon}{2M_2}, \sqrt{\frac{\varepsilon}{L_2}}\right\}, \frac{M_2}{L_2}\sqrt{\frac{1}{6n}}\right\},$$

в итоге можно получить

$$\delta \leq \frac{\varepsilon^{3/2}}{16R\sqrt{L_2 n}}.$$

**Замечание 3 (техника двойного сглаживания).** Сопоставляя условия, при которых оценки таблицы 2 справедливы для обратной связи вида (8) и вида (9), несложно заметить, что в случае (8) в теореме 5 требовалось выполнение дополнительного условия гладкости (12). Замечание 2 хорошо поясняет почему при $\delta > 0$ это условие может понадобиться. Однако уже довольно давно (Б.Т. Поляк, 1983 [21]; см. также [19]) было подмечено, что и без условия (12) можно использовать безградиентные двухточечные методы. Для этого вместо (8) формируется

$$\nabla_x f\left(x^k, \xi^k\right) = \frac{n}{\tau_2}\left(f\left(x^k + \tau_1\tilde{e}_1^k + \tau_2 e_2^k, \xi^k\right) - f\left(x^k + \tau_2\tilde{e}_1^k, \xi^k\right)\right)e_2^k,$$

где $\tilde{e}_1^k \in RB_2^n(1)$, т.е. $\tilde{e}_1^k$ – случайный вектор, равномерно распределенный на евклидовом шаре единичного радиуса в пространстве $\mathbb{R}^n$, $e_2^k \in RS_2^n(1)$ и $\left\{\tilde{e}_1^k, e_2^k, \xi^k\right\}_k$ – независимы в совокупности. В работе удалось получить условие на уровень шума, при котором оценки из второго столбца таблицы 2, отвечающего негладкому случаю, для обратной связи (8) остаются по-прежнему верными, т.е. верными без условия (12),

$$\delta \le \frac{\varepsilon^2}{56 M_2 R n^{3/2}}.$$

При этом

$$\tau_1 = \frac{\varepsilon}{4M_2},\ \tau_2 = \frac{\varepsilon}{4M_2 n}.$$

Насколько нам известно, ранее не был известны теоретически обоснованные способы выбора $\tau_1$, $\tau_2$. Как и следовало ожидать, из-за отсутствия гладкости, приведенное здесь условие на допустимый уровень шума $\delta$ получилось более жестким, чем в гладком случае – см. замечание 2.

**Пример 3.** Пусть на воображаемой планете живут “оптимизаторы”. Каждый из них стремиться жить с наименьшими потерями нервов. Планета характеризуется выпуклой $M_2$-липшицевой функцией $f(x)$, аргумента $x \in Q \subseteq \mathbb{R}^n$, где $Q$ – выпуклое множество, являющееся множеством допустимых стратегий оптимизаторов. Считаем, что этому множеству соответствует выбор $p = 2$. Функция $f(x)$ обозначает нервные потери оптимизатора за день, если оптимизатор в этот день использовал стратегию $x \in Q$. Однако значение функции $f(x)$ не наблюдаемы. Оптимизатор, использовавший в $k$-й день стратегию $x^k$, может наблюдать лишь зашумленную реализацию значения этой функции $f(x^k, \xi^k)$.

Рассматривается два способа жить на такой планете: 1) жить одному, 2) жить вдвоем. В случае 1) оптимизатор может получать каждый день только обратную связь вида (7), а вот в случае 2), за счет сговора, каждому можно рассчитывать на обратную связь вида (8). Из теоремы 5 и замечания 3 (в онлайн варианте) получаем, что в случае 1) необходимо (и достаточно) прожить $N = C_1 B^2 M_2^2 R^2 n^2 / \varepsilon^4$ дней, чтобы обеспечить выполнение условия

$$E[\mathrm{Regret}_N] = \frac{1}{N}\sum_{k=1}^{N} E\left[f(x^k)\right] - \min_{x\in Q}\frac{1}{N}\sum_{k=1}^{N} f(x) \le \varepsilon.$$

В случае 2) для этой же цели достаточно прожить $N = C_2 M_2^2 R^2 n / \varepsilon^2$ дней. Здесь $C_1, C_2 \le 100$. Интересно заметить, что если жить втроем, вчетвером и т.д., то это лишь немого улучшает оценку $N = C_2 M_2^2 R^2 n / \varepsilon^2$. А именно, если $2k$ оптимизаторов живут вместе и могут согласовывать свои стратегии, исходя из общедоступной информации, то $N \approx (C_2/k) M_2^2 R^2 n / \varepsilon^2$. □

## 4. *Стохастический зеркальный спуск и нащупывание равновесия в популяционной игре загрузки*

В данном разделе, базируясь на работе [22], рассматривается приложение результатов раздела 2 к решению задач поиска равновесного распределения потоков в транспортных сетях.

Рассматривается транспортная сеть, которую можно представить ориентированным графом $\langle V, E \rangle$, где $V$ – множество вершин, а $E$ – множество ребер. Обозначим множество пар источник-сток через $OD \subseteq V \otimes V$; $d_w$ – корреспонденция, отвечающая паре $w$; $x_p$ – поток по пути $p$; $P_w$ – множество путей, отвечающих корреспонденции $w$, $P = \bigcup_{w \in OD} P_w$ – множество всех путей. Обозначим через $H$ – максимальное число ребер в пути из $P$. Затраты на прохождения ребра $e \in E$ описываются неубывающей и ограниченной в рассматриваемом диапазоне значений функцией

$$0 \le \tau_e(f_e) \le \tilde{M},$$

где $f_e$ – поток по ребру $e$,

$$f_e(x) = \sum_{p \in P} \delta_{ep} x_p, \; \delta_{ep} = \begin{cases} 1, & e \in p \\ 0, & e \notin p \end{cases}.$$

Положим $\breve{M} = \tilde{M} \cdot H$. Введем затраты на проезд по пути $p$

$$G_p(x) = \sum_{e \in E} \tau_e(f_e(x)) \delta_{ep}.$$

Введем также множество

$$X = \left\{ x \ge 0 : \sum_{p \in P_w} x_p = d_w, w \in OD \right\},$$

и функцию, порождающее потенциальное векторное поле $G(x)$:

$$\Psi(x)=\sum_{e\in E}\int_0^{f_e(x)}\tau_e(z)dz.$$

Основное свойство этой функции заключается в том, что

$$\nabla\Psi(x)=G(x).$$

Будем считать, что число пользователей транспортной сети большое

$$d_w:=d_w\cdot\bar{N},\ \bar{N}\gg 1,\ w\in OD,$$

но в функциях затрат это учитывается

$$\tau_e(f_e):=\tau_e\left(f_e/\bar{N}\right).$$

Таким образом, далее под $d_w$, $x$, $f$ будем понимать соответствующие прошкалированные по $\bar{N}$ величины [23].

Выберем корреспонденцию $w\in OD$ и рассмотрим пользователя транспортной сети, соответствующего этой корреспонденции. Стратегией пользователя является выбор одного из возможных путей следования $p\in P_w$. Будем считать, что пользователь мало что знает об устройстве транспортной системы и о формировании своих затрат. Все что доступно пользователю на шаге $k+1$ – это история затрат на разных путях, соответствующих его корреспонденции, на всех предыдущих шагах $\left\{l_p^r=\left\{G_p\left(x^r\right)\right\}_{p\in P_w}\right\}_{r=1}^k$. Для простоты рассуждений мы не зашумляем эту информацию, считая что доступны точные значения имевших место затрат. Все последующие рассуждения можно обобщить и на случай зашумленных данных.

Допуская, что $0\le\left\{l_p^k\right\}\le\breve{M}$ могут выбираться враждебно, пользователь стремиться действовать оптимальным образом, то есть так, как предписывает стратегия из примера 2 с $i=p$, $n=n_w=|P_w|$. Заметим, что при некоторых дополнительных оговорках (см. [9]) случайный выбор пути согласно рандомизированному варианту МЗС из примера 2 может быть осуществлен за время $\mathrm{O}(|E|)$, что не зависит от $n$, которое может быть намного больше, как это имеет место, например, для манхетенской сети.

Представим себе, что остальные пользователи ведут себя аналогичным образом, но независимо (в вероятностном плане) друг от друга. Тогда в пределе $\bar{N} \to \infty$ такая стохастическая марковская динамика в дискретном времени вырождается в детерминированную динамику $x^k \in \prod_{w \in OD} S_{n_w}(1)$, $k = 1, ..., N$ в дискретном времени, и имеет место следующий результат.

**Теорема 6.** *Пусть*

$$\bar{x}^N = \frac{1}{N}\sum_{k=1}^{N} x^k \text{ , } \Psi_* = \Psi(x_*), x_* \in \operatorname{Arg}\min_{x \in X} \Psi(x).$$

*Тогда*

$$\Psi(\bar{x}^N) - \Psi_* \le \frac{M}{\sqrt{N}} \frac{\max_{w \in OD}\{\ln n_w\}}{\sqrt{2\min_{w \in OD}\{\ln n_w\}}} \left( \sum_{w \in OD} d_w^2 + 1 \right).$$

Решение $x_*$ часто называют равновесием Нэша(–Вардропа) в описанной выше классической модели равновесного распределения потоков по путям (Beckmann–McGuire–Winsten, 1954 [24]). Такая "онлайн" интерпретация этого равновесия, насколько нам известно, приводится впервые.

К аналогичному равновесию можно было прийти и по-другому. Второй подход воспроизводит logit-динамику в соответствующей описанной модели популяционной игре загрузки. Приводимый далее способ рассуждения пополняет аналогичные известные ранее результаты, полученные сразу в непрерывном времени [23].

Пусть теперь $l$-й пользователь, принадлежащий корреспонденции $w \in W$, независимо от всех остальных пользователей на шаге $k+1$ с вероятностью $1 - \lambda/N$ выбирает путь $p^{l,k}$, который использовал на шаге $k = 0, ..., tN$, а с вероятностью $\lambda/N$ ($\lambda > 0$) решает изменить путь и выбирает (возможно, тот же самый) зашумленный кратчайший путь

$$p^{l,k+1} = \arg\max_{q \in P_w} \left\{ -G_q(x^k) + \xi_q^{l,k+1} \right\},$$

$$G_q(x) = \sum_{e \in E} \tau_e(f_e(x))\delta_{eq} \text{ , } f_e(x) = \sum_{p \in P} \delta_{ep} x_p \text{, } \delta_{ep} = \begin{cases} 1, & e \in p \\ 0, & e \notin p \end{cases};$$

где независимые случайные величины $\xi_q^{l,k+1}$ одинаково распределены

$$P\left(\xi_q^{l,k+1} < \zeta\right) = \exp\left\{-e^{-\zeta/\gamma - E}\right\}, \ \gamma > 0.$$

Такое распределение обычно называю распределение Гумбеля или двойным экспоненциальным распределением. Возникает оно здесь не случайно, Andersen–de Palma–Thisse, 1992 [25]. В математической экономике имеется целое направление Discrete choice theory, которое объясняет использование именно такого распределения. А именно, распределение Гумбеля можно объяснить исходя из идемпотентного аналога центральной предельной теоремы (вместо суммы случайных величин – максимум) для независимых случайных величин с экспоненциальным и более быстро убывающим правым хвостом. Распределение Гумбеля возникает в данном контексте, например, если при принятии решения водитель собирает информацию с большого числа разных (независимых) зашумленных источников, ориентируясь на худшие прогнозы по каждому из путей. Также это распределение удобно тем, что есть явная формула

$$P\left(p^{l,k+1}=p \middle| \text{агент решил "поменять" путь}\right)=\frac{\exp\left(-G_p\left(x^k\right)/\gamma\right)}{\sum_{q\in P_w}\exp\left(-G_q\left(x^k\right)/\gamma\right)}.$$

Если взять $E\approx 0.5772$ – константа Эйлера, то

$$E\left[\xi_q^{l,k+1}\right]=0,\ D\left[\xi_q^{l,k+1}\right]=\gamma^2\pi^2/6.$$

Пусть $N\to\infty$. Тогда описанная марковская динамика в дискретном времени перейдет в марковскую динамику в непрерывном времени. Если $t\to\infty$, то полученный таким образом марковский процесс в виду эргодичности выйдет на стационарную меру (теорема типа Санова)

$$\sim\exp\left(-\frac{\bar{N}}{\gamma}\cdot\left(\Psi_\gamma\left(x\right)+o_{\bar{N}}\left(1\right)\right)\right),\ \bar{N}\gg 1,$$

где

$$\Psi_\gamma\left(x\right)=\Psi\left(x\right)+\gamma\sum_{w\in W}\sum_{p\in P_w}x_p\ln x_p\ .$$

При $\bar{N}\to\infty$ эта мера экспоненциально быстро концентрируется в окрестности «стохастического» равновесия [26], которое определяется из решения задачи

$$\Psi_\gamma\left(x\right)\to\min_{x\in X}. \tag{13}$$

Поскольку $\Psi_\gamma\left(x\right)$ – строго выпуклая функция, то стохастическое равновесие всегда единственно.

Отметим, что при $\gamma\to 0+$ описанная выше logit-динамика вырождается в одну из самых популярных и естественных динамик в популяционной теории игр (W. Sandholm, 2010 [23]) – динамику наилучших ответов. При этом стохастическое равновесие вырождается в равновесие Нэша(–Вардропа). Если последних равновесий много, то описанный выше предельный переход $\gamma\to 0+$ можно понимать как естественный способ отбора единственного (Е.В. Гасникова, 2012 [27]). “Естественный”, потому что связан с естественно интерпретируемой энтропийной регуляризацией.

Приведенная выше трактовка равновесия с точки зрении популяционной теории игр допускает и эволюционную (дарвиновскую) интерпретацию. Если под разными типами корреспонденций понимать разные типы популяций. Под пользователями, использующими внутри заданной корреспонденции–популяции заданный маршрут, понимать заданный вид внутри заданной популяции, то динамика наилучших ответов – это просто естественный отбор (борьба за существование), в результате которого “выживет” такая конфигурация видов, которая максимально приспособлена. Таким образом, равновесие Нэша(–Вардропа) может быть тут также проинтерпретировано, как результат естественного отбора, т.е. как равновесие Дарвина. Причем, поскольку $x_* \in \operatorname{Arg}\min_{x \in X} \Psi(x)$, то имеет место эволюционная оптимальность с функционалом $\Psi(x)$ (В.Н. Разжевайкин, 2012 [28]).

К сожалению, хорошо интерпретируемый способ поиска $x_*$ из теоремы 6 не является наиболее подходящим численным методом поиска $x_*$. Тому как именно стоит численно искать $x_*$ планируется посвятить отдельные работы. В завершении данного раздела приведены необходимые и также интерпретируемые выкладки по переформулировки исходной задачи (13), с целью возможности ее последующего эффективного численного решения.

Запишем, следуя [29] двойственную задачу к задаче (13), используя обозначение $\operatorname{dom} \sigma^*$ – область определения сопряженной к $\sigma$ функции:

$$\min_{f,x}\left\{\sum_{e \in E}\sigma_e\left(f_e\right)+\gamma\sum_{w \in OD}\sum_{p \in P_w}x_p\ln\left(x_p/d_w\right):\ f=\Theta x,\ x \in X\right\}=$$

$$=\min_{f,x}\left\{\sum_{e \in E}\max_{t_e \in \operatorname{dom}\sigma_e^*}\left[f_e t_e-\sigma_e^*\left(t_e\right)\right]+\gamma\sum_{w \in OD}\sum_{p \in P_w}x_p\ln\left(x_p/d_w\right):\ f=\Theta x,\ x \in X\right\}=$$

$$=\max_{t \in \operatorname{dom}\sigma^*}\left\{\min_{f,x}\left[\sum_{e \in E}f_e t_e+\gamma\sum_{w \in OD}\sum_{p \in P_w}x_p\ln\left(x_p/d_w\right):\ f=\Theta x,\ x \in X\right]-\sum_{e \in E}\sigma_e^*\left(t_e\right)\right\}=$$

$$=-\min_{t \in \operatorname{dom}\sigma^*}\left\{\gamma\psi\left(t/\gamma\right)+\sum_{e \in E}\sigma_e^*\left(t_e\right)\right\}, \qquad (14)$$

где

$$\psi\left(t\right)=\sum_{w \in OD}d_w\psi_w\left(t\right),\ \psi_w\left(t\right)=\ln\left(\sum_{p \in P_w}\exp\left(-\sum_{e \in E}\delta_{ep}t_e\right)\right),$$

$$f=-\nabla\gamma\psi(t/\gamma),\ x_p=d_w\frac{\exp\left(-\frac{1}{\gamma}\sum_{e\in E}\delta_{ep}t_e\right)}{\sum_{q\in P_w}\exp\left(-\frac{1}{\gamma}\sum_{e\in E}\delta_{eq}t_e\right)},\ p\in P_w, \tag{15}$$

для

$$\tau_e(f_e)=\bar{t}_e\cdot\left(1+\rho\cdot\left(\frac{f_e}{\bar{f}_e}\right)^{\frac{1}{\mu}}\right),$$

$$\sigma_e^*(t_e)=\sup_{f_e\ge 0}\left((t_e-\bar{t}_e)\cdot f_e-\bar{t}_e\cdot\frac{\mu}{1+\mu}\cdot\rho\cdot\frac{f_e^{1+\frac{1}{\mu}}}{\bar{f}_e^{\frac{1}{\mu}}}\right)=\bar{f}_e\cdot\left(\frac{t_e-\bar{t}_e}{\bar{t}_e\cdot\rho}\right)^{\mu}\frac{(t_e-\bar{t}_e)}{1+\mu}.$$

В приложениях наиболее популярны BPR-функции, в которых $\mu=1/4$.

Собственно, формула (14) есть не что иное, как отражение формулы $f=-\nabla\gamma\psi(t/\gamma)$ и связи $t_e=\tau_e(f_e)$, $e\in E$. Действительно, по формуле Демьянова–Данскина–Рубинова

$$\frac{d\sigma_e^*(t_e)}{dt_e}=\frac{d}{dt_e}\max_{f_e\ge 0}\left\{t_ef_e-\int_0^{f_e}\tau_e(z)dz\right\}=f_e:\ t_e=\tau_e(f_e).$$

В свою очередь, формула $f=-\nabla\gamma\psi(t/\gamma)$ может интерпретироваться, как следствие соотношений $f=\Theta x$ и формулы распределения Гиббса (15) (logit-распределения)

$$x_p=d_w\frac{\exp\left(-\frac{1}{\gamma}\sum_{e\in E}\delta_{ep}t_e\right)}{\sum_{q\in P_w}\exp\left(-\frac{1}{\gamma}\sum_{e\in E}\delta_{eq}t_e\right)},\ p\in P_w.$$

При такой интерпретации связь задачи (14) с logit-динамикой, порождающей стохастические равновесия, наиболее наглядна. Действительно,

$$f\xrightarrow{t_e=\tau_e(f_e)}\underbrace{t\xrightarrow{(15)}x\xrightarrow{f=\Theta x}f}_{t\xrightarrow{f=-\nabla\gamma\psi(t/\gamma)}f}.$$

Поиск неподвижной точки в этой цепочке, как раз и сводится к решению задачи (14).

Работа поддержана грантом Президента РФ МД-1320.2018.1 и РФФИ 18-29-03071 мк.

## *Литература*

1. *Немировский А.С., Юдин Д.Б.* Сложность задач и эффективность методов оптимизации. – М.: Наука, 1979. – 384 с.
2. *Nemirovski A.* Lectures on modern convex optimization analysis, algorithms, and engineering applications. Philadelphia: SIAM, 2016. – URL: http://www2.isye.gatech.edu/~nemirovs/Lect_ModConvOpt.pdf
3. *Juditsky A., Lan G., Nemirovski A., Shapiro A.* Stochastic approximation approach to stochastic programming // SIAM Journal on Optimization. – 2009. – V. 19. – № 4. – P. 1574–1609.
4. *Двуреченский П.Е., Гасников А.В., Лагуновская А.А.* Параллельные алгоритмы и оценки вероятностей больших уклонений в задачах стохастической выпуклой оптимизации // Сиб. ЖВМ. – 2018. – Т. 21. – № 1. – С. 47–53.
5. *Nesterov Y., Vial J.-Ph.* Confidence level solution for stochastic programming // Automatica. – 2008. – V. 44. – no. 6. – P. 1559–1568.
6. *Hazan E.* Introduction to online convex optimization // Foundations and Trends® in Optimization. 2016. V. 2: No. 3-4. P. 157–325.
7. *Juditsky A., Nemirovski A.* First order methods for nonsmooth convex large-scale optimization, I, II. // Optimization for Machine Learning. // Eds. S. Sra, S. Nowozin, S. Wright. – MIT Press, 2012.
8. *Гасников А.В., Крымова Е.А., Лагуновская А.А., Усманова И.Н., Федоренко Ф.А.* Стохастическая онлайн оптимизация. Одноточечные и двухточечные нелинейные многорукие бандиты. Выпуклый и сильно выпуклый случаи // Автоматика и Телемеханика. – 2017. – № 2. – С. 36–49.
9. *Lugosi G., Cesa-Bianchi N.* Prediction, learning and games. – New York: Cambridge University Press, 2006.

10. *Hazan E., Kale S.* Beyond the regret minimization barrier: Optimal algorithms for stochastic strongly-convex optimization // JMLR. – 2014. – V. 15. – P. 2489–2512.
11. *Nesterov Yu.* Primal-dual subgradient methods for convex problems // Math. Program. Ser. B. – 2009. – V. – 120(1). – P. 261–283.
12. *Гасников А.В., Лагуновская А.А., Усманова И.Н., Федоренко Ф.А.* Безградиентные прокс-методы с неточным оракулом для негладких задач выпуклой стохастической оптимизации на симплексе // Автоматика и телемеханика. – 2016. – № 10. – С. 57–77.
13. *Горбунов Э.А., Воронцова Е.А., Гасников А.В.* О верхней оценке математического ожидания нормы равномерно распределённого на сфере вектора и явлении концентрации равномерной меры на сфере // Матем. зам. – 2019. (в печати) https://arxiv.org/pdf/1804.03722.pdf
14. *Ledoux M.* Concentration of measure phenomenon. – Providence, RI, Amer. Math. Soc., 2001 (Math. Surveys Monogr. V. 89).
15. *Баяндина А.С., Гасников А.В., Лагуновская А.А.* Безградиентные двухточечные методы решения задач стохастической негладкой выпуклой оптимизации при наличии малых шумов не случайной природы // Автоматика и телемеханика. – 2018. – №. 8. – С. 38–49.
16. *Agarwal A., Dekel O., Xiao L.* Optimal algorithms for online convex optimization with multi-point bandit feedback // Proceedings of 23 Annual Conference on Learning Theory. –2010. – P. 28–40.
17. *Bubeck S., Eldan R.* Multi-scale exploration of convex functions and bandit convex optimization // e-print, 2015. – URL: http://research.microsoft.com/en-us/um/people/sebubeck/ConvexBandits.pdf
18. *Bubeck S., Cesa-Bianchi N.* Regret analysis of stochastic and nonstochastic multi-armed bandit problems // Foundation and Trends in Machine Learning. – 2012. – V. 5. – № 1. – P. 1–122.

19. *Duchi J.C., Jordan M.I., Wainwright M.J., Wibisono A.* Optimal rates for zero-order convex optimization: the power of two function evaluations // IEEE Transaction of Information. – 2015. – V. 61. – № 5. – P. 2788–2806.
20. *Shamir O.* An optimal algorithm for bandit and zero-order convex optimization with two-point feedback // Journal of Machine Learning Research. – 2017. – V. 18. – P. 1–11.
21. *Поляк Б.Т.* Введение в оптимизацию. – М.: Наука, 1983. – 384 с.
22. *Гасников А.В., Лагуновская А.А., Морозова Л.Э.* О связи имитационной логит динамики в популяционной теории игр и метода зеркального спуска в онлайн оптимизации на примере задачи выбора кратчайшего маршрута // Труды МФТИ. – 2015. – Т. 7. – № 4. – С. 104–113.
23. *Sandholm W.* Population games and Evolutionary dynamics. Economic Learning and Social Evolution. – MIT Press; Cambridge, 2010.
24. *Beckmann M., McGuire C.B., Winsten C.B.* Studies in the economics of transportation. – RM-1488. Santa Monica: RAND Corporation, 1955.
25. *Andersen S.P., de Palma A., Thisse J.-F.* Discrete choice theory of product differentiation. – MIT Press; Cambridge, 1992.
26. *Sheffi Y.* Urban transportation networks: Equilibrium analysis with mathematical programming methods. – N.J.: Prentice–Hall Inc., Englewood Cliffs, 1985.
27. *Гасникова Е.В.* Моделирование динамики макросистем на основе концепции равновесия: дис. канд. физ.-мат. наук: 05.13.18. – М.: МФТИ, 2012. – 90 с.
28. *Разжевайкин В.Н.* Анализ моделей динамики популяций. Учебное пособие. – М.: МФТИ, 2010. – 174 с.
29. *Гасников А.В., Гасникова Е.В., Мациевский С.В., Усик И.В.* О связи моделей дискретного выбора с разномасштабными по времени популяционными играми загрузок // Труды МФТИ. – 2015. – Т. 7. – № 4. – С. 129–142.